\documentclass[notitlepage,leqno,10pt]{article}

\usepackage[applemac]{inputenc}
\usepackage{amsfonts}
\usepackage{amsmath}
\usepackage{amssymb}
\usepackage{amsrefs}

\def\neweq#1{\begin{equation}\label{#1}}
\def\endeq{\end{equation}}
\def\eq#1{(\ref{#1})}

\usepackage{color}
\definecolor{green2}{rgb}{0,0.6,0}

\newcommand{\R}{\mathbb{R}}

\newcommand{\proof}[1]{\par\smallskip\noindent{\bf Proof#1.}}
\newcommand{\qed}{\penalty 500\hfill$\square$\par\medskip}

\newcommand{\eps}{\varepsilon}

\newcommand{\super}[1]{\overline{#1} }

\def\XXint#1#2#3{{\setbox0=\hbox{$#1{#2#3}{\int}$}
     \vcenter{\hbox{$#2#3$}}\kern-.5\wd0}}

\newtheorem{lem}         {Lemma}

\newtheorem{thm}    [lem]{Theorem}
\newtheorem{rem}    [lem]{Remark}

\title{Regularity of the extremal solutions for the Liouville system}
\date{\today}
\author{L. Dupaigne$^{1,2}$, A. Farina$^{1}$ and B. Sirakov$^{3}$}
\begin{document}
\maketitle

{\begin{center}
$^{1}${\small
LAMFA, UMR CNRS 7352, Universit\'e Picardie Jules Verne \\33, rue St Leu, 80039 Amiens, France \\
\smallskip
$^2$corresponding author, {\tt louis.dupaigne@math.cnrs.fr}\\
\smallskip
$^3${\small Pontif\`{i}cia Universidade Cat\'olica do Rio de Janeiro (PUC-Rio)
Departamento de Matem\'atica
Rua Marqu\^{e}s de S\~{a}o Vicente, 225, G\'avea
Rio de Janeiro - RJ,
CEP 22451-900, Brasil\\}}
\end{center}}

In this short note, we study the smoothness of the extremal solutions to the following system of equations:
\begin{equation} \label{liouville} 
\left\{
\begin{aligned}
-\Delta u&= \mu e^v&\qquad\text{in $\Omega$,}\\
-\Delta v&= \lambda e^{u}&\qquad\text{in $\Omega$,}\\
u&=v=0&\qquad\text{on $\partial\Omega$,}
\end{aligned}
\right.
\end{equation}
where $\lambda,\mu>0$ are parameters and $\Omega$ is a smoothly bounded domain of $\R^N$, $N\ge1$.
As shown by M. Montenegro (see \cite{marcelo}), there exists a limiting curve $\Upsilon$ in the first quadrant of the $(\lambda,\mu)$-plane serving as borderline for existence of classical solutions of \eqref{liouville}. He also proved the existence of a weak solution $u^*$ for every $(\lambda^*,\mu^*)$ on the curve $\Upsilon$ and left open the question of its regularity. Following standard terminology (see e.g. the books \cite{dupaigne}, \cite{egg} for an introduction to this vast subject), $u^*$ is called an extremal solution. 
Our result is the following.
\begin{thm}\label{th main} Let $1\le N\le 9$. Then, extremal solutions to \eqref{liouville} are smooth. 
\end{thm}
\begin{rem}
C. Cowan (\cite{cowan}) recently obtained the same result under the further assumption that  $(N-2)/8<\lambda/\mu<8/(N-2)$. 
\end{rem}

Any extremal solution $u^*$ is obtained as the increasing pointwise limit of a sequence of regular solutions $(u_{n})$ associated to parameters $(\lambda_{n},\mu_{n})=(1-1/n)(\lambda^*,\mu^*)$.  In addition, see \cite{marcelo},  $u_{n}$ is stable in the sense that the principal eigenvalue of the linearized operator associated to \eqref{liouville}  is nonnegative. In other words, there exist $\lambda_{1}\ge0$ and two positive functions $\varphi_{1},\psi_{1}\in C^2(\super\Omega)$ such that
\begin{equation}\label{linearized} 
\left\{
\begin{aligned}
-\Delta \varphi_{1}-g'(v)\psi_{1}&= \lambda_{1}\varphi_{1}&\qquad\text{in $\Omega$,}\\
-\Delta \psi_{1}-f'(u)\varphi_{1}&= \lambda_{1}\psi_{1}&\qquad\text{in $\Omega$.}\\
\varphi_{1}&=\psi_{1}=0&\qquad\text{on $\partial\Omega$,}
\end{aligned}
\right.
\end{equation}
where, in the context of \eqref{liouville},  $g(v)=e^{v}$ and $f(u)=e^{u}$. This motivates the following useful inequality.

Let $f,g$ denote two nondecreasing $C^1$ functions and consider the more general system
\begin{equation} \label{main} 
\left\{
\begin{aligned}
-\Delta u&= g(v)&\qquad\text{in $\Omega$,}\\
-\Delta v&= f(u)&\qquad\text{in $\Omega$,}\\
u&=v=0&\qquad\text{on $\partial\Omega$.}
\end{aligned}
\right.
\end{equation}
 
\begin{lem}\label{lemma 1} Let $ N\ge 1$ and let $(u,v)\in C^2(\super\Omega)^2$ denote a stable solution of \eqref{main}. Then, for all $\varphi\in C^1_{c}(\Omega)$, there holds
\begin{equation} \label{stability}
\int_{\Omega}\sqrt{f'(u)g'(v)}{\varphi^2}\;dx\le \int_{\Omega}\left| \nabla \varphi \right|^2\;dx 
\end{equation}   
\end{lem}

\begin{rem}As we just learnt, the same inequality has been obtained independently by  C. Cowan and N. Ghoussoub. See \cite{cg}.
\end{rem}

\proof{}
Since $(u,v)$ is stable, there exist $\lambda_{1}\ge0$ and two positive functions $\varphi_{1},\psi_{1}\in C^2(\super\Omega)$ solving \eqref{linearized}. 
Given $\varphi\in C^1_{c}(\Omega)$, multiply the first equation in \eqref{linearized} by $\varphi^2/\varphi_{1}$ and integrate. Then,
\neweq{sys1}
\begin{aligned}
\int_\Omega g'(v) \frac{\psi_1}{\varphi_1}\varphi^2\;dx&\le\int_\Omega \frac{\varphi^2}{\varphi_1}(-\Delta \varphi_1)\\
&=-\int_\Omega |\nabla \varphi_1|^2\Big(\frac{\varphi}{\varphi_1}  \Big)^2+2\int_\Omega \frac{\varphi}{\psi_1}\nabla \varphi\nabla \varphi_1\\
&=-\int_\Omega \Big|\frac{\varphi}{\varphi_1} \nabla \varphi_1-\nabla \varphi\Big|^2+\int_\Omega |\nabla \varphi|^2\leq \int_\Omega |\nabla \varphi|^2.
\end{aligned}
\endeq
Working similarly with the second equation, we also have
\begin{equation} \label{stab2}
\int_{\Omega}f'(u)\frac{\varphi_{1}}{\psi_{1}}\varphi^2\;dx\le \int_{\Omega}\left| \nabla \varphi\right|^2\;dx 
\end{equation}
\eqref{stability} then follows by combining the Cauchy-Schwarz inequality and \eqref{sys1}- \eqref{stab2}.    
\hfill\qed

Thanks to the inequality \eqref{stability}, we obtain the following estimate. 
\begin{lem} \label{lemma 2} Let $ N\ge 1$. There exists a universal constant $C>0$  such that any stable solution of \eqref{liouville} satisfies
\begin{equation}\label{l1}  
\int e^{u+v}\;dx \le C \left| \Omega \right| \left(\frac\lambda\mu+\frac\mu\lambda\right).
\end{equation} 
\end{lem}

\proof{}
Multiply the second equation in \eqref{liouville} by $e^{v}-1$ and integrate.
\begin{align}\label{a1} 
\lambda\int_{\Omega}e^{u+v}dx \ge
\lambda\int_{\Omega}e^{u}(e^v-1)dx&=
\int_{\Omega}\nabla v \nabla (e^{v}-1)\;dx\\\nonumber
&=4\int_{\Omega}\left| \nabla(e^{v/2}-1)\right|^2 \;dx.
\end{align} 
Using \eqref{stability} with test function $\varphi=e^{v/2}-1$, it follows that
\begin{align} \label{a2}
\lambda\int_{\Omega}e^{u+v}\; dx &\ge 4\sqrt{\lambda\mu}\int_{\Omega}e^{\frac{u+v}2}(e^{v/2}-1)^2\;dx\\ \nonumber
&\ge  4\sqrt{\lambda\mu}\int_{\Omega}e^{\frac{u+v}2}e^{v}\;dx -  8\sqrt{\lambda\mu}\int_{\Omega}e^{\frac{u+v}2}e^{v/2}\;dx.
\end{align}  
By Young's inequality, $e^{v/2}=\frac1{\sqrt2} e^{v/2}\cdot \sqrt 2\le \frac14 e^v +1$. So,
$$
\int_{\Omega}e^{\frac{u+v}2}e^{v/2}\;dx\le \frac14 \int_{\Omega}e^{\frac{u+v}2}e^{v}\;dx + \int_{\Omega}e^{\frac{u+v}2}\;dx.
$$
Plugging this in \eqref{a2}, we obtain 
\begin{equation} \label{a3}
\lambda\int_{\Omega}e^{u+v}\; dx+ 8\sqrt{\lambda\mu}\int_{\Omega}e^{\frac{u+v}2}\;dx\ge 2 \sqrt{\lambda\mu}\int_{\Omega}e^{\frac{u+v}2}e^{v}\;dx.
\end{equation}  
Similarly,
\begin{equation} \label{a4}
\mu\int_{\Omega}e^{u+v}\; dx+ 8\sqrt{\lambda\mu}\int_{\Omega}e^{\frac{u+v}2}\;dx\ge 2 \sqrt{\lambda\mu}\int_{\Omega}e^{\frac{u+v}2}e^{u}\;dx.
\end{equation} 
Multiply \eqref{a3} and \eqref{a4} to get
\begin{multline}\label{a5}  
\lambda\mu\left(\int_{\Omega}e^{u+v}\; dx\right)^2 +64\lambda\mu\left(\int_{\Omega}e^{\frac{u+v}2}\;dx\right)^2 + 8\sqrt{\lambda\mu}(\lambda+\mu)
\int_{\Omega}e^{u+v}\; dx\int_{\Omega}e^{\frac{u+v}2}\;dx
\ge\\4\lambda\mu\int_{\Omega}e^{\frac{u+v}2}e^{u}\;dx\int_{\Omega}e^{\frac{u+v}2}e^{v}\;dx.
\end{multline} 
Using Young's inequality, the left-hand side in the above inequality is bounded above by
\begin{equation} \label{a6} 
2\lambda\mu\left(\int_{\Omega}e^{u+v}\; dx\right)^2 + C(\lambda+\mu)^2\left(\int_{\Omega}e^{\frac{u+v}2}\;dx\right)^2,
\end{equation} 
where $C$ is a universal constant. In addition, 
by the Cauchy-Schwarz inequality, 
\begin{equation} \label{a7} 
\int_{\Omega}e^{\frac{u+v}2}e^{u}\;dx\int_{\Omega}e^{\frac{u+v}2}e^{v}\;dx \ge 
\left(\int_{\Omega}e^{u+v}\; dx\right)^2.
\end{equation} 
Plugging \eqref{a7} in \eqref{a6} and remembering that \eqref{a6} is an upper bound   of the left-hand side in \eqref{a5}, we obtain 
\begin{equation} \label{a8}
C(\lambda+\mu)^2\left(\int_{\Omega}e^{\frac{u+v}2}\;dx\right)^2\ge 
2\lambda\mu\int_{\Omega}e^{\frac{u+v}2}e^{u}\;dx\int_{\Omega}e^{\frac{u+v}2}e^{v}\;dx.
\end{equation}  
By the Cauchy-Schwarz inequality and \eqref{a7}, we have
\begin{align} \label{a9} 
\left(\int_{\Omega}e^{\frac{u+v}2}\;dx\right)^2&\le \left| \Omega \right| \int_{\Omega}e^{u+v}\;dx\\\nonumber
&\le \left| \Omega \right|\left(\int_{\Omega}e^{\frac{u+v}2}e^{u}\;dx\int_{\Omega}e^{\frac{u+v}2}e^{v}\;dx\right)^{1/2}.
\end{align}   
Using \eqref{a9} in \eqref{a8}, we obtain
\begin{equation} 
C\frac{(\lambda+\mu)^2}{\lambda\mu}\left| \Omega \right|\ge  \left(\int_{\Omega}e^{\frac{u+v}2}e^{u}\;dx\int_{\Omega}e^{\frac{u+v}2}e^{v}\;dx\right)^{1/2}.
\end{equation}  
Applying once more \eqref{a7}, we obtain the desired estimate.  
\hfill\qed
We can now prove Theorem \ref{th main}.

\

\noindent {\bf Step 1. Case $1\le N\le3$. }
It is enough to treat the case $N=3$, the cases $N=1,2$ being easier. 
By \eqref{a1} and \eqref{l1}, $e^{v/2}-1$ is bounded in $H^1_{0}(\Omega)$ (with a uniform bound with respect to $ \lambda $ and $\mu$). By the Sobolev embedding, it follows that $e^v$ is bounded in $L^{\frac N{N-2}}(\Omega)$. By \eqref{a1} and elliptic regularity, $u$ is bounded in $W^{2, \frac N{N-2}}$.  For $N=3$, $\frac N{N-2}>\frac N2$. By Sobolev's embedding, we deduce that $u$ is bounded, and so must be $v$.  This implies the desired conclusion for the corresponding extremal solution. 


\noindent {\bf Step 2. General case. }
We adapt a method introduced in \cite{dggw}. Fix $\alpha>1/2$ and multiply the first equation in \eq{liouville} by $e^{\alpha u}-1$. Integrating over $\Omega$, we obtain
$$
\mu\int_{\Omega} \left( e^{\alpha u}-1\right)e^v\;dx = \alpha \int_{\Omega}e^{\alpha u}\vert\nabla u\vert^2\;dx=
\frac4\alpha \int_{\Omega}\left\vert\nabla\left(e^{\frac{\alpha u}2}-1\right)\right\vert^2\;dx
$$
By \eqref{stability},
$$
\sqrt{\lambda\mu}\int_{\Omega}e^{\frac{u+v}2}\left(e^{\frac{\alpha u}{2}}-1\right)^2\;dx\le \int_{\Omega}\left\vert\nabla\left(e^{\frac{\alpha u}2}-1\right)\right\vert^2\;dx.
$$ 
Combining these two inequalities, we deduce that
\begin{equation} \label{2161} 
\sqrt{\lambda\mu}\int_{\Omega}e^{\frac{u+v}2}\left(e^{\frac{\alpha u}{2}}-1\right)^2\;dx \le \frac{\alpha}{4}\mu\int_{\Omega} \left( e^{\alpha u}-1\right)e^v\;dx\end{equation} 
Hence,
\begin{equation} \label{2161b} 
\sqrt{\lambda\mu}\int_{\Omega}e^{\frac{2\alpha+1}2u}e^{\frac v2}\;dx \le \frac{\alpha}{4}\mu\int_{\Omega} e^{\alpha u}e^v\;dx + 2\sqrt{\lambda\mu}\int_{\Omega}e^{\frac{\alpha+1}{2}u}e^{\frac v2}\;dx
\end{equation} 
Let us estimate the terms on the right-hand side. By H\"older's inequality,
\begin{equation}\label{rhs1} 
\int_{\Omega}e^{\alpha u}e^{v}\;dx \le \left(\int_{\Omega} e^{\frac{2\alpha+1}2u}e^{\frac v2}\;dx\right)^{\frac{2\alpha-1}{2\alpha}}\left(\int_{\Omega} e^{\frac u2}e^{\frac{2\alpha+1}2 v}\;dx\right)^{\frac1{2\alpha}}
\end{equation}
Given $\eps>0$, it also follows from Young's inequality that
$$\int_{\Omega}e^{\frac{\alpha+1}{2}u}e^{\frac v2}\;dx\le \frac\eps{2} \sqrt \frac\mu{\lambda}\int_{\Omega}e^{\alpha u}e^{v}\;dx + \frac1{2\eps} \sqrt \frac\lambda{\mu}\int_{\Omega}e^{u}\;dx.
$$
Using \eqref{l1}, we deduce that 
\begin{equation} \label{rhs2} 
\int_{\Omega}e^{\frac{\alpha+1}{2}u}e^{\frac v2}\;dx\le\frac\eps{2} \sqrt \frac\mu{\lambda}\int_{\Omega}e^{\alpha u}e^{v}\;dx +\frac1{2\eps} \sqrt \frac\lambda{\mu}C \left| \Omega \right| \left(\frac\lambda\mu+\frac\mu\lambda\right).\end{equation}
where $C$ is the universal constant of Lemma \ref{lemma 2}.

So, gathering \eqref{2161b}, \eqref{rhs1}, \eqref{rhs2}, and letting 
$$
X=\int_{\Omega}e^{\frac{2\alpha+1}2u}e^{\frac v2}\;dx\quad\text{and}\quad Y=\int_{\Omega}e^{\frac{2\alpha+1}2v}e^{\frac u2}\;dx,
$$
we obtain
$$
\sqrt{\lambda\mu }\, X \le \left(\frac{\alpha}{4}+\eps\right)\mu\, X^{\frac{2\alpha-1}{2\alpha}}Y^{\frac1{2\alpha}} + C \frac\lambda\eps \left| \Omega \right| \left(\frac\lambda\mu+\frac\mu\lambda\right).$$


By symmetry, we also have

$$
\sqrt{\lambda\mu }\, Y \le \left(\frac{\alpha}{4}+\eps\right)\lambda\, Y^{\frac{2\alpha-1}{2\alpha}}X^{\frac1{2\alpha}} + C \frac\mu\eps \left| \Omega \right| \left(\frac\lambda\mu+\frac\mu\lambda\right).$$


Multiplying these inequalities, we deduce that
$$\left(1-\left(\frac{\alpha}{4}+\eps\right)^2\right)X\,Y \le C_1  \left(\frac\lambda\mu+\frac\mu\lambda\right)^2
\left(1+ X^{\frac{2\alpha-1}{2\alpha}}Y^{\frac1{2\alpha}} + Y^{\frac{2\alpha-1}{2\alpha}}X^{\frac1{2\alpha}} \right).
$$
where $C_1 =  C \frac{\left| \Omega \right|} \eps \left(\frac{\alpha}{4}+\eps\right)>0$. 
Hence, for every $\alpha<4$, either $X$ or $Y$ must be bounded (with a uniform bound with respect to $ \lambda $ and $\mu$).

\noindent Without loss of generality, $\lambda\ge\mu$ and by the maximum principle, $v\ge u$. It follows that $e^{u}$ is bounded in $L^p(\Omega)$ for every $p=\alpha+1<5$. Using standard elliptic regularity, the result follows.\hfill\qed
\bibliographystyle{amsalpha}

\end{document}